\documentclass[letterpaper, 10 pt, conference]{IEEEconf}
\IEEEoverridecommandlockouts

\usepackage{cite}
\usepackage{amsmath,amssymb,amsfonts,commath}
\usepackage{algorithmic}
\usepackage{graphicx}
\usepackage{textcomp}
\usepackage{xcolor}
\usepackage{url}
\usepackage{siunitx}
\usepackage{makecell}
\graphicspath{{figures/}}
\usepackage{graphicx}
\usepackage[]{subfigure}
\usepackage{bm}
\usepackage{hyperref}
\usepackage{multirow}
\usepackage{booktabs}
\usepackage{amsmath}
\usepackage{mathtools}  
\usepackage{xfrac} 

\def\BibTeX{{\rm B\kern-.05em{\sc i\kern-.025em b}\kern-.08em
    T\kern-.1667em\lower.7ex\hbox{E}\kern-.125emX}}

\newcommand{\ez}{\ensuremath{\bm{\mathbf{e}}_{z}}}

\newcommand{\er}{\ensuremath{\bm{\mathbf{e}}_{r}}}
\newcommand{\ephi}{\ensuremath{\bm{\mathbf{e}}_{\varphi}}}
\newcommand{\ebeta}{\ensuremath{\bm{\mathbf{e}}_{\beta}}}

\newcommand{\eX}{\ensuremath{\bm{\mathbf{e}}_{X}}}
\newcommand{\eY}{\ensuremath{\bm{\mathbf{e}}_{Y}}}
\newcommand{\eZ}{\ensuremath{\bm{\mathbf{e}}_{Z}}}

\newcommand{\F}{\ensuremath{\bm{\mathbf{F}}}}

\newcommand{\cbeta}{\ensuremath{{\,\cos\beta\,}}}
\newcommand{\tbeta}{\ensuremath{{\,\tan\beta\,}}}

\newcommand{\calpha}{\ensuremath{{\,\cos\alpha\,}}}
\newcommand{\ctheta}{\ensuremath{{\,\cos\theta\,}}}
\newcommand{\cgamma}{\ensuremath{{\,\cos\gamma\,}}}
\newcommand{\sbeta}{\ensuremath{{\,\sin\beta\,}}}

\newcommand{\salpha}{\ensuremath{{\,\sin\alpha\,}}}
\newcommand{\stheta}{\ensuremath{{\,\sin\theta\,}}}
\newcommand{\sgamma}{\ensuremath{{\,\sin\gamma\,}}}

\begin{document}

\bibliographystyle{IEEEtran}

\title{Automatic Circular Take-off and Landing of Tethered Motorized Aircraft

\thanks{This research is supported by  FCT/MCTES(PIDDAC), through projects 2022.02320.PTDC-KEFCODE, and 2022.02801.PTDC-UPWIND-ATOL (https://doi.org/10.54499/2022.02801.PTDC).}
}

\author{
Sérgio Vinha,
Gabriel M. Fernandes,
Huu Thien Nguyen, \\
Manuel C.R.M. Fernandes, and 
Fernando A.C.C. Fontes.
\thanks{Sérgio Vinha, Gabriel M. Fernandes, Huu Thien Nguyen, Manuel C.R.M. Fernandes, and Fernando A.C.C. Fontes are with SYSTEC-ISR-ARISE and the Department of Electrical and Computer Engineering, Faculty of Engineering, University of Porto, Rua Dr. Roberto Frias, 4200-465 Porto, Portugal. Emails: 
        {\tt\small svinha@fe.up.pt},
        {\tt\small gabriel.fernandes@upwind.pt},
        {\tt\small nguyen@fe.up.pt},
        {\tt\small mcrmf@fe.up.pt},
        {\tt\small faf@fe.up.pt}
        }%
}
\maketitle
\begin{abstract}

We consider a motorized aircraft tethered to a central anchorage point in a configuration similar to a control line model airplane. 
For this system, we address the problem of automatic take-off and landing (ATOL) with a circular path, {whose} center and radius are defined by the anchorage point and the tether length, respectively.
We propose a hierarchical control architecture for ATOL and discuss the controllers designed for each control layer and for each of the flight phases.
Simulation results are reported, showing the viability of the approach, but also showing the limitations on the maximum altitude attainable with a fixed-tether length. 
The tethered aircraft and the proposed ATOL control architecture are to be used in an Airborne Wind Energy System.

\end{abstract}

\begin{keywords}
Aircraft Control, 
Tethered Airplanes,
Automatic Take-off and Landing, Hierarchical Control,
Airborne Wind Energy.
\end{keywords}

\section{Introduction}

This article addresses the development of an automatic take-off and landing framework for a tethered airplane. 
In particular, we develop a control system for circular take-off and landing of a self-propelled, fixed-wing, tethered aircraft.
The main application that we foresee for the proposed technology, is to include it in an Airborne Wind Energy System.

Airborne Wind Energy Systems (AWES) are devices that convert wind energy into electricity using autonomous aircraft attached to the ground by a tether \cite{loyd_crosswind_1980, schmehl_airborne_2018}. These devices can harvest wind energy at high altitudes, where the wind is stronger and more consistent, being able to generate electricity from a yet unexplored renewable energy resource.
One of the main challenges in the development of AWES into a commercially viable and competitive renewable energy technology is the ability to operate safely, reliably, and autonomously for long periods of time in several weather and environmental conditions (see e.g., the European Commission report \cite{directorate-general_for_research_and_innovation_european_commission_study_2018} and a recent survey \cite{vermillion_electricity_2021}). Most existing technology demonstrators still rely on supervised operation, especially in the take-off and landing phases, and are, therefore, not fully autonomous. To achieve fully autonomous operation it is crucial to develop reliable Automatic Take-Off and Landing (ATOL) schemes for tethered aircraft.

Some AWES technologies use fixed-wing aircraft, while others rely on soft wings. Depending on the characteristics of these wings,  different automatic take-off and landing (ATOL) approaches can be adopted \cite{vermillion_electricity_2021}. Among several AWES using fixed-wing, tethered and motorized aircraft, we can distinguish three methods: linear, vertical and circular TOL techniques \cite{fagiano_take-off_2017}. The least investigated technique of the three, and the one studied in this article, is the Circular Take-off and Landing (CTOL). The CTOL scheme has been developed in two different variants. One of the variants, the most investigated one, uses a rotating platform or rigid arm to generate enough airspeed on the kite \cite{sieberling_discussion_2011, zanon_rotational_2013, geebelen_experimental_2013, rieck_comparison_2017}. Other variants use kites equipped with landing gear and propellers to produce the circular motion \cite{cherubini_preliminary_2017, cherubini_advances_2017, fernandes_take-off_2023} during the take-off and landing phases, similarly to the kite studied here. However, a detailed analysis of the dynamical system and the control system is still to be done. It also is curious to note that Miles Loyd, one of the pioneers of Airborne Wind Energy,  describes in its seminal patent \cite{loyd_wind_1981} an apparatus based on a circular launch using a vehicle coupled to a circular rail system.

In this paper, we propose a control system for an automatic circular take-off and landing that can be adopted by a tethered self-motorized aircraft. 
We use a hierarchical control architecture. In the top layer, we design a supervisory
controller that is responsible for governing the transition
between flight phases, for path-planning, and for setting the references to the lower-level controllers at each phase of operation. 
The controllers designed for each phase range from simple PID, designed for one control-variable, to multivariable optimal regulators for the locally linearized systems.
The developed framework has been tested in simulations and in a small-scale prototype.
The results show the viability of the approach to take-off and attain a certain altitude as well as to the process of landing. The results also show the limits on the maximum altitude attainable with a fixed-tether length as a function of the path radius.

\begin{table}[!t]
\centering
\caption{Nomenclature} \label{tab:nomenclature}
\begin{tabular}{@{}ll}
{$A$} & {wing area $\left[\rm{m}^2\right]$}
\\
{$b$} & {wingspan $\left[\rm{m}\right]$}
\\
{$c_{\rm{D}},c_{\rm{L}}$} & {aerodynamic drag and lift coefficients}
\\
{$\F_{\rm{D}},F_{\rm{D}}$} & {aerodynamic drag vector and magnitude $\left[\rm{N}\right]$}
\\
{$\F_{\rm{L}},F_{\rm{L}}$} & {aerodynamic lift vector and magnitude $\left[\rm{N}\right]$}
\\
{$\F_{\rm{p}},F_{\rm{p}}$} & {propeller thrust vector and magnitude $\left[\rm{N}\right]$}
\\
{$\F_{\rm{i}},F_{\rm{i}}$} & {inertial force vector and magnitude in a rotating frame $\left[\rm{N}\right]$}
\\
{$\F_{\rm{t}},F_{\rm{t}}$} & {tether force vector and magnitude $\left[\rm{N}\right]$}
\\
{$g$} & {gravitational acceleration $\left[\rm{m\,s^{-2}}\right]$}
\\
{$m$} & {mass $\left[\rm{kg}\right]$}
\\
{$\rho$} & {air density $\left[\rm{kg\,m^{-3}}\right]$}
\\
{$h$} & {height $\left[\rm{m}\right]$}
\\
{$\mathbf{p}$}  & {kite position} 
\\
{$\mathbf{V_a},V_a$ } & {kite airspeed vector and magnitude $\left[\rm{m\,s^{-1}}\right]$}
\\
{$V_w$ } & {wind velocity $\left[\rm{m\,s^{-1}}\right]$}
\\
{$\mathbf{X}, \mathbf{U}$} & {state and control vectors}
\\
{$r,\varphi,\beta$} & {spherical coordinates $\left[\rm{m}\right]$,$\left[^\circ,\rm{rad}\right]$,$\left[^\circ,\rm{rad}\right]$}
\\
{$\alpha,\gamma$} & {angle-of-attack and flight-path angle $\left[^\circ,\rm{rad}\right]$}
\\
{$\phi,\theta,\psi$} & {roll, pitch and yaw angles $\left[^\circ,\rm{rad}\right]$}
\\
{$\tau$} & {tangent plane at the kite position}
\\
{$\omega_p,\omega_q,\omega_r$} & {roll, pitch and yaw rates $\left[\rm{^\circ\,s^{-1}},\rm{rad\,s^{-1}}\right]$}

\\

\end{tabular}
\end{table}
\section{Tethered Aircraft Model}

The position and motion of a tethered aircraft  can be conveniently represented in spherical coordinates 
$(r,\varphi,\beta)$, centered at the tether anchorage point, where $r$ is the radial distance, $\varphi$ is the azimuth angle, and $\beta$ the elevation angle (from the horizontal plane), having coordinate basis  $(\er,\ephi,\ebeta)$. See nomenclature in Table \ref{tab:nomenclature}.

The position, velocity, and acceleration vectors are \cite{riley2006mathematical}:
\begin{align}
&\mathbf{p}
= 
\begin{bmatrix}
r  \\
\varphi\\
\beta
\end{bmatrix}
, 
\quad
\mathbf{\dot{p}} =
\begin{bmatrix}
\dot{r} \\
r \cbeta \dot{\varphi}  \\
r \dot{\beta}
\end{bmatrix}
, 
\mbox{\ and} \nonumber\\
&
\mathbf{\ddot{p}} =
\begin{bmatrix}
\ddot{r} \\
r\cbeta\ddot{\varphi} \\
r\ddot{\beta}
\end{bmatrix}
+
\begin{bmatrix}
- r \dot{\beta}^2 - r \dot{\varphi}^2 \cos^2\beta \\
2\dot{r}\dot{\varphi}\cbeta - 2r\dot{\varphi}\dot{\beta}\sbeta \\
2\dot{r}\dot{\beta} + r\dot{\varphi}^2\cbeta\sbeta
\end{bmatrix}
,
\end{align}
where the last term in $\mathbf{\ddot{p}}$ is due to the inertial forces $\F_{\rm{i}}$ arising from the use of a rotating frame and equals $- \F_{\rm{i}}/m$.

If we assume for the moment that the tether is nonelastic, always taut, and with constant length, the kite is constrained to a 2D motion on a surface of a sphere of radius $r$, equal to the tether length. 
In this case, it is convenient to consider a plane $\mathbf{\tau}$ that is tangential to the spherical surface at the kite position, which is the span of the basis vectors $\ephi,\ebeta$.

\begin{figure}[!t]
\centering
\subfigure{\includegraphics[scale=0.33, trim={10cm 5cm 5.5cm 3cm}, clip]{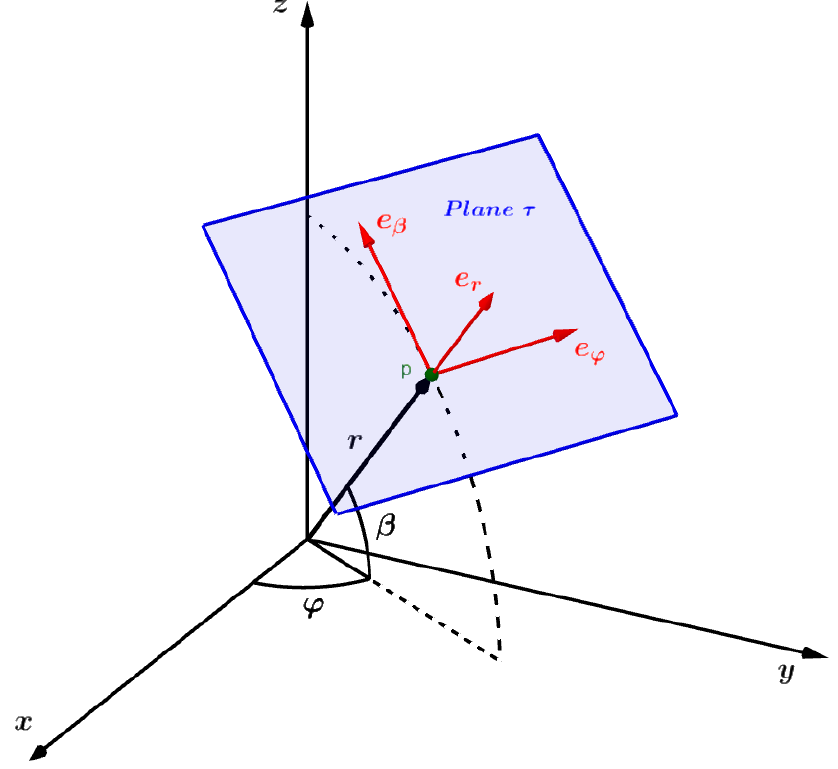}}
\qquad 
\subfigure{\includegraphics[width=.17\textwidth, trim={0cm -1cm 0cm 0cm}, clip]{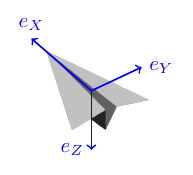}}
\caption{Spherical and Body coordinate systems.}
\label{fig:coordinate_system}
\end{figure}

We also consider a coordinate frame attached to the kite body, with basis $(\eX,\eY,\eZ)$, where $\eX$ is the aircraft's longitudinal axis pointing to its nose, 
$\eY$ is the transversal axis pointing out of the right wing, and $\eZ$ is the aircraft's vertical axis pointing down from its belly (see Fig. \ref{fig:coordinate_system}).

The pitch angle $\theta$  (defined to be the angle of $\eX$ with the horizontal plane), the flight-path angle $\gamma$ (defined to be the angle of the airspeed vector $\mathbf{V_a}$ with the horizontal plane),
and the aircraft angle-of-attack $\alpha$ (defined to be the angle of $\eX$ with $\mathbf{V_a}$) are related by 
$
\theta=\gamma+\alpha.
$

\begin{figure}[h!]
    \centering
    \includegraphics[width=.3\textwidth]{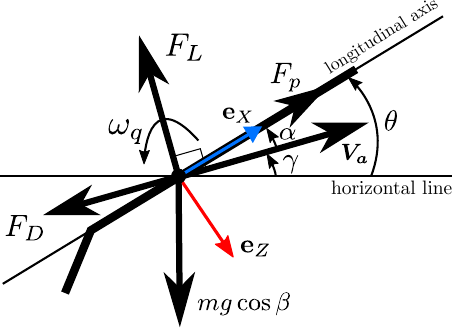}
    \caption{Airplane longitudinal model for circular tethered flight.}
    \label{fig:kite}
\end{figure}

\begin{figure}[!t]
\centering
\includegraphics[trim={3cm 3.5cm 2cm 5cm},clip,width=\columnwidth]{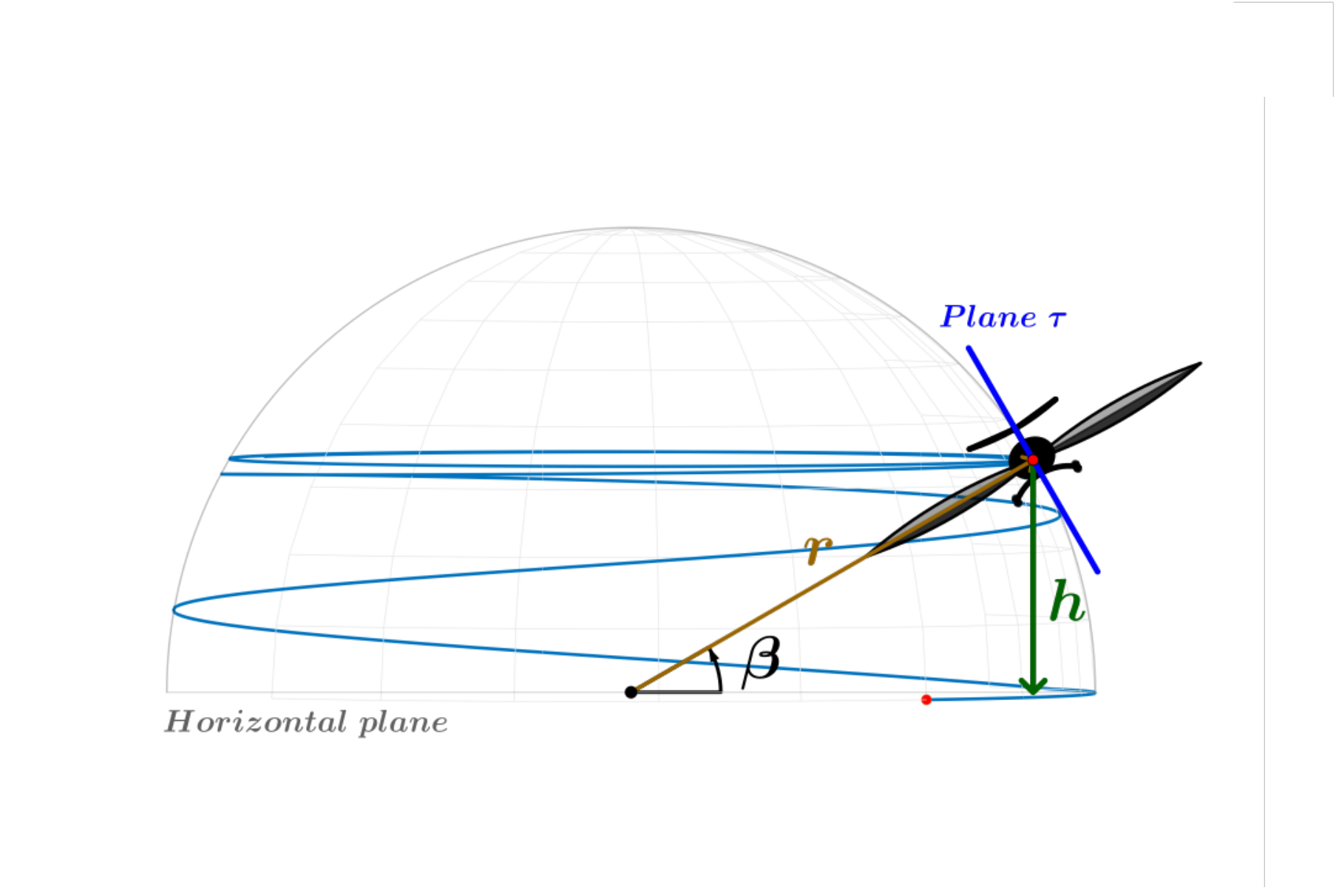}
\caption{Circular trajectory during take-off (in blue).}
\label{fig:sphere}
\end{figure}

The kinematics of the aircraft satisfy 
\begin{equation}
\begin{aligned} \label{eq:kin}
r \cbeta \dot{\varphi} =& V_a \cgamma \\
r \dot{\beta} =& V_a \sgamma \\
\dot{\theta} =& \:\omega_q,
\end{aligned}
\end{equation}
where $\omega_q$ is the pitch rate.
Taking the derivative of the kinematics \eqref{eq:kin}, we have
\begin{equation}
\begin{aligned} \label{eq:acc}
r \ddot{\beta}  = &  \sgamma{\dot{V} }_a +V_a \,\cgamma\dot{\gamma} \\
r \cbeta \ddot{\varphi}  = & \cgamma{\dot{V} }_a - V_a\sgamma\dot{\gamma} \\
& + \frac{{V_a }^2}{r}\tan\beta \cgamma \sgamma. \\
\end{aligned}
\end{equation}
The forces acting on the aircraft comprise the propeller thrust $\F_p=F_p \: \eX$, the weight $\F_g$ with magnitude $mg$ pointing down,  the aerodynamic lift $ \F_L=- F_L \: \ez$, the aerodynamic drag $\F_D= - F_D \: \eX$, the tether pull $\F_t= - F_t \: \er$, and also the inertial forces $\F_i$ when a rotating frame is considered, yielding
 $
   \F= [F_r \, F_{\varphi} \, F_{\beta}] ^T=\F_L + \F_D +\F_p  + \F_g +  \F_t + \F_i. 
 $
We consider that the tether point of attachment to the aircraft is at the wing tip aligned with the center of mass (see Fig. \ref{fig:sphere}) and if we assume a no-wind situation, then, the axis $\eY$ aligns with $\er$, the velocity vector  $\mathbf{V}_a$ and the aerodynamic forces are in the tangential plane $\mathbf{\tau}$. In this setting,
$F_D =\frac{1}{2}\rho A  {V_a} ^2 c_{\rm{D}}{(\alpha)} $ and  $F_L =\frac{1}{2}\rho A {V_a} ^2 c_{\rm{L}}{(\alpha)}$.  

Equating $ m \mathbf{\ddot{p}}$ to all external forces and projecting into
 the tangential plane $\mathbf{\tau}$,  with $\dot{r}=0, \ddot{r}=0$, we have in the 
$\ephi$, $\ebeta$ directions
\begin{align} \label{eq:dyn}
  m r \cbeta \ddot{\varphi} &= F_{\varphi} \nonumber\\
  m r \dot{\beta} &= F_{\beta},
\end{align}
with 
\begin{align}
   F_{\varphi} =& -F_L\sgamma - F_D\cgamma + F_p\ctheta + 2mr\dot{\varphi}\dot{\beta}\sbeta \nonumber \\ 
  F_{\beta} =&  F_L\cgamma - F_D\sgamma + F_p\stheta - mg\cbeta \\
  &- mr\dot{\varphi}^2\cbeta\sbeta.   \nonumber
\end{align}
Combining these last equations with \eqref{eq:acc} and solving for  $(\dot{V}_a,\dot{\gamma})$, we obtain
\begin{equation}
\begin{aligned} \label{eq:model}
m \dot{V}_a =& -{F_D} + {F_p}\,\calpha -mg\cbeta \sgamma \\
m V_a \dot{\gamma} =& {F_L}+{F_p}\salpha -mg\cbeta \cgamma \\
 & -\frac{m {V_a}^2 }{r}\,\tan\beta  \cgamma,
\end{aligned}
\end{equation}
which together with the kinematics \eqref{eq:kin} is the analogous to circular tethered flight of the well-known
   longitudinal model of a fixed-wing plane (see \cite{nguyen_trajectory_2023, beard_small_2012} and Fig. \ref{fig:kite}).

We  assume that the thrust force $F_p$ and the pitch rate $\omega_q$ are variables that we can manipulate directly (in fact, via the autopilot). 
Then, the state and control vectors considered are, respectively,
\begin{equation}
\mathbf{X} = \begin{bmatrix}
\varphi, 
\beta, 
V_a, 
\gamma, 
\theta
\end{bmatrix}^\top, \quad
\mathbf{U}
= 
\begin{bmatrix}
F_p,  
\omega_q
\end{bmatrix}^\top,
\end{equation}
and the state-space model is  given by
\begin{equation} \label{statespace}
\mathbf{\dot{X}} = f(\mathbf{X},\mathbf{U}), 
\end{equation}
where $f$ is taken from equations \eqref{eq:kin} and \eqref{eq:model}.

\section{Control Architecture}
We propose a control architecture for a Circular Take-off and Landing (CTOL) of a rigid-wing tethered aircraft, envisaging application to Airborne Wind Energy Systems (AWES). We consider a kite equipped with a landing gear and propeller, only to be used during the Take-off and Landing (TOL) phases.

Parting from rest, the kite accelerates on the ground while performing a circular trajectory with a fixed radius around the point of attachment to the ground station. 
After reaching a set speed for take-off, the aircraft increases its pitch angle and elevates until it reaches a desired altitude. 
From that point onward, the kite is able to ascend, having the tether reeling out (thus increasing the radius of the circular trajectory) and positioning itself downwind at an appropriate point to begin its power productive phases.
The landing procedure should be similar but in a reversed order. Parting from its power production motion, the kite will start following a horizontal circular trajectory centered at the ground station while reeling in the tether to the initial radius and descending to a predefined altitude suitable for landing.
After landing, it will reduce its speed until it comes to a full stop on the ground.

Below, in Subsection \ref{sub:SC}, we introduce the high-level Supervisory Controller; in Subsection \ref{sec:Controllers}, we describe the lower-level Airspeed and Altitude Controllers.

\subsection{Supervisory Control}\label{sub:SC}
In order to design a controller for automatic take-off and landing, it is convenient to divide the full operation of the system into several phases with different control requirements, references, and methods. Therefore, we design a  Supervisory Controller that is responsible for governing the transition between phases, for path-planning, and for setting the references to the lower-level controllers at each phase of operation. 

\begin{figure}[!t]
    \centering
    \includegraphics[width=0.48\textwidth]{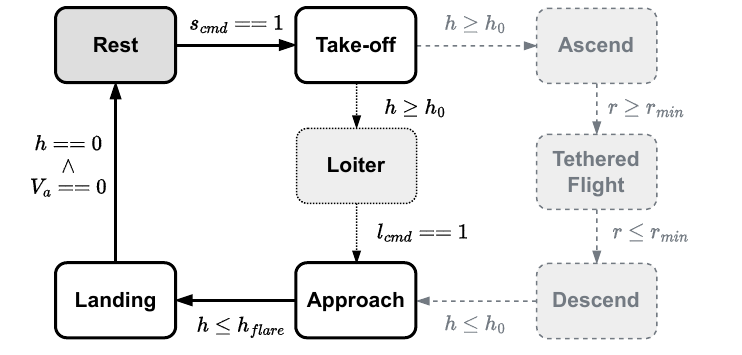}
    \caption{Overall Control Strategy of an AWE system.}
    \label{fig:SC}
\end{figure}

The diagram in Fig. \ref{fig:SC} represents the different phases of operation and the conditions in which the system transitions from one to the next.
The goal of this paper is to develop a controller for the TOL phases, thus we will only discuss results and simulations that refer to the unblurred phases in Fig. \ref{fig:SC}. The controllers for the Ascend, Tethered Flight, and Descend phases (portrayed in blurred and dashed lines) will be discussed elsewhere.
To link all the phases during the simulation, we add a Loiter phase connecting the Take-off and Approach phases where the kite flies in circles with a constant altitude and tether length (portrayed in pointed lines)
Also, in this work, we focus on the Kite Module control and do not discuss the Ground-station Module control 
(for the ground station controller, we point to \cite{uppal_cascade_2021, arshad_uppal_ground_2022}).

\begin{figure}[!t]
    \centering
    \includegraphics[width=0.48\textwidth,trim={0.5cm 0cm 0.5cm 0cm},clip]{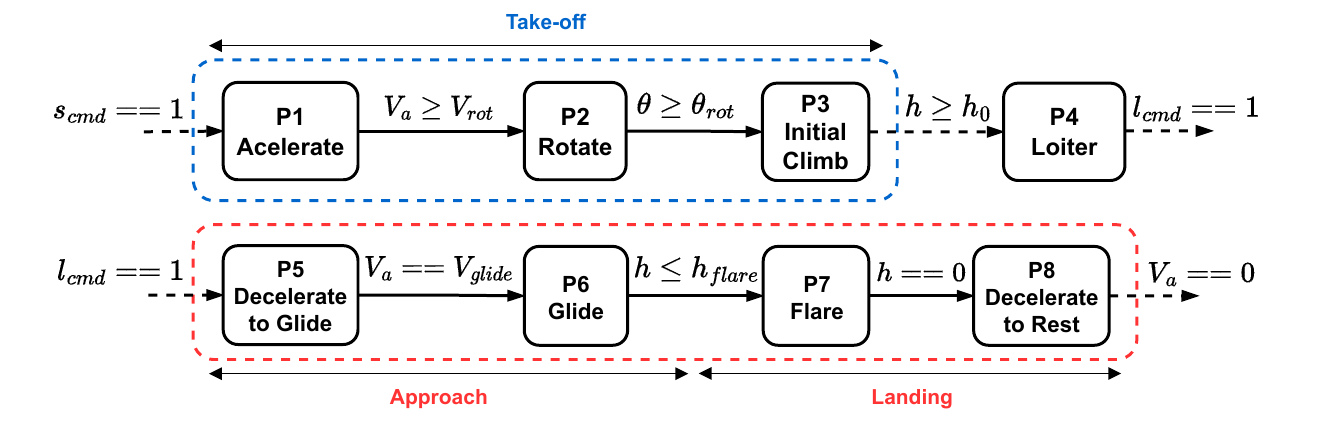}
    \caption{The Take-off, Loiter, Approach, and Landing phases with its sub-phases.}
    \label{fig:takeoff_landing}
\end{figure}

When in Rest, the system will await for an operator command and/or suitable wind and weather conditions for energy production and will jump to the following phase of Take-off.
This phase comprises three sub-phases, depicted in Fig. \ref{fig:takeoff_landing} with the blue dashed rectangle. 
These sub-phases move sequentially the kite from a motionless position on the ground to a level flight, ready to ascend and later begin generating electricity.
Firstly, in the sub-phase of Acceleration, the kite will increase its speed while on the ground until it reaches a predefined speed  $v_{rot}$. 
After reaching $v_{rot}$, the kite will pass onto the sub-phase Rotation, in which it will tilt upwards until it reaches a reference pitch angle $\theta_{rot}$. This will confer to the kite its maximum lift and move it to the Initial Climb phase where it will maintain its speed and pitch angle while increasing its altitude.

When a threshold altitude $h_0$ is reached, the system switches to the Ascend phase. At this stage, the kite increases its altitude, increases the tether length,  and directs its circular trajectory downwind, thus defining a wider circular and quasi-crosswind motion. When the tether length reaches the minimum value needed for the pumping cycle operation, $r_{min}$, the system will jump to the Tethered Flight phase and start generating electrical power.

Whenever there is the need to retrieve the kite (\textit{e.g.} for maintenance purposes or due to hazardous weather conditions) the system initiates its Descend phase. Here the tether is reeled back in and the circular trajectory will return to a horizontal plane at a reduced altitude until it reaches $h_0$ and switches to the Approach phase.

During the \textit{Approach} mode, the aircraft flies in a circular and near horizontal path with a fixed tether length. It starts to decelerate to a certain velocity $V_{glide}$, while maintaining the altitude. Then, it starts to glide and the altitude decreases until a certain minimum altitude $h_{flare}$.
Finally, \textit{Landing} phase is reached. Here, in a flare maneuver, the angle-of-attack of the aircraft increases and the speed decreases right before touchdown. On the ground, the speed continues to decrease until the aircraft stops at a resting position.
These two last phases are highlighted in Fig. \ref{fig:takeoff_landing} with a red dashed rectangle.

\subsection{Altitude and Airspeed Controllers}
\label{sec:Controllers}

\begin{figure}[!t]
    \centering
    \includegraphics[width=0.48\textwidth,trim={0.5cm 0cm 0.5cm 0cm},clip]{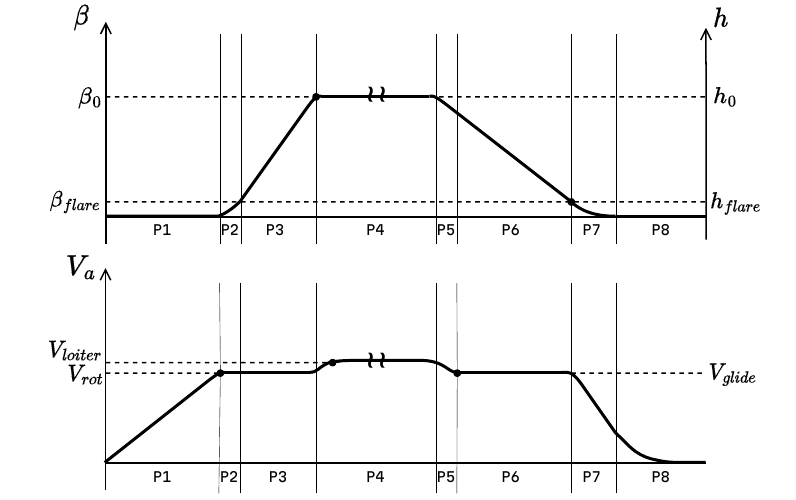}
    \caption{Altitude and Airspeed reference profile during all phases.}
    \label{fig:phases}
\end{figure}

The main objective of the ATOL control strategy is to guarantee that the kite airspeed, altitude and attitude values have the desired behavior according to predefined reference values at each phase.
To achieve such an objective, we use two different control methodologies. 
In the phases where we are mainly concerned with changing rapidly the value of one variable -- as in P1, P2, P5, and P7 -- we use  simple PID controllers, which perform competently. 
In phase P8, since we do not have any active breaking system, we simply set the propulsion to zero and wait for the drag and friction to stop the aircraft.
In phases P3, P4, and P6, we want to maintain several state components near their steady state, and therefore we use a multivariable controller, a Linear Quadratic Regulator (LQR), computed for the system linearized around a steady midpoint of the trajectory in each phase (see Table \ref{table:target_states} and Fig. \ref{fig:phases}). The PID and LQR controllers are further detailed in the next subsections. 

In the different flight phases, two main angles-of-attack (AoA) are considered: 
the steady flight AoA $\alpha_0$, in which the ratio $c_{\rm{L}}/ c_{\rm{D}}$ is maximized, and the maximum lift AoA $\alpha_{\rm{L}}$ for which $c_{\rm{L}}$ is maximized. While the steady flight AoA $\alpha_0$ is used in phase P4-Loiter (and also in tethered flight energy producing phases), the maximum lift AoA $\alpha_{\rm{L}}$ is used in phases P3 and P6. During phase P2, P5 and P7, the AoA will vary. During P2-Rotate we aim for a pitch angle $\theta_{rot}=\alpha_{\rm{L}}-\alpha_0$ and the speed $V_{rot}$ is selected in such a way that the lift with AoA $\alpha_0$ is not enough to overcome the aircraft weight, but the lift with AoA $\alpha_{\rm{L}}$ can elevate the aircraft (i.e.  $F_{\rm{L}}(\alpha_0,V_{rot}) < mg < F_{\rm{L}}(\alpha_{\rm{L}},V_{rot})$).

\begin{table}[!t]
\caption{Target states and respective controller for each phase} 
\label{table:target_states}
\centering
\begin{tabular}{llc}
\textbf{Phase} - \textbf{Name}  & \textbf{Target State} (conditions)   & \textbf{Controller} \\
\toprule
\textbf{P1} - Accelerate     &  $\begin{aligned}&V_a  = V_{rot}, \theta = 0 \\
&(\alpha=\alpha_0,\gamma=0, \beta = 0)\end{aligned}$ &  2 PID    \\ 
\midrule

\textbf{P2} - Rotate         &   $\begin{aligned}
& V_a  = V_{rot}, \theta=\theta_{rot} \\
&(\dot{\gamma},\dot{\alpha} \geq 0, \dot{\beta} \geq =0)
\end{aligned}$  &    2 PID\\
\midrule

\textbf{P3} - Initial Climb  &    $\begin{aligned}&V_a = V_{climb}, \gamma = \gamma_{climb}\\ &\alpha = \alpha_{\rm{L}}  \end{aligned}$                     & LQR  \\    
\midrule

\textbf{P4} - Loiter   &          $\begin{aligned}&\beta=\beta_0, V_a = V_{loiter}\\&\gamma=0, \theta = 0, \alpha=\alpha_0  \end{aligned}$               &    LQR \\
\midrule

\textbf{P5} - $\begin{aligned}&\text{Deccelerate} \\&\text{to Glide} \end{aligned}$     &  $\begin{aligned}&V_a  = V_{glide}, \gamma  = 0 \\& (\dot{\alpha} \geq 0)\end{aligned}$ & 2 PID  \\          
\midrule

\textbf{P6} - Glide          &  $\begin{aligned}&V_a = V_{glide}, \gamma = \gamma_{glide} \\&\alpha = \alpha_{\rm{L}} \end{aligned}$ & LQR\\
\midrule

\textbf{P7} - Flare          &    $\theta=\theta_{flare}$ ($F_p=0$)                  &     PID \\
\midrule

\textbf{P8} - $\begin{aligned}&\text{Deccelerate} \\&\text{to Rest} \end{aligned}$   &   $V_a=0$ ($F_p=0$)              &        --  \\
\bottomrule

\end{tabular}
\end{table}

\subsubsection{PID Controllers}

The PID controllers are used to control the speed $V_a$ by actuating on the thrust $F_p$ or to control the attitude (measuring $\theta$ or $\gamma$) by actuating on the pitch rate $\omega_q$. 
In phases P1, P2 and P5, we use two PID controllers with references (target values) given in Table \ref{table:target_states}. In phase P7-Flare, we use just one PID to control the attitude; the velocity is naturally decreased by making the thrust $F_p$ equal to zero and by the augmented drag due to the higher AoA.
The tuning of the PID gains was done as in \cite{fernandes_take-off_2023} and the values for our simulation can be seen in Table \ref{table:controller_gains}.

\subsubsection{Linear-Quadratic Regulator}
The state considered in this section does not include the azimuth angle $\varphi$.
In circular take-off and landing, we have an endless runway and there is no need to control the azimuth position $\varphi$. Moreover, the evolution of $\beta$, $V_a$, $ \gamma$, and $\theta$ described by the last four component equalities in (8) do not depend on $\varphi$. 
Therefore, we can omit $\varphi$ from the state-space considered for control purposes; we will just integrate the first equality to draw the trajectories as in Fig. \ref{fig:sphere}. The state and control considered is then
$
\mathbf{x} = ( \beta, V_{a}, \gamma, \theta)
$, 
$
\mathbf{u} = ( F_{p}, w_{q} ),
$
which satisfy $\dot{\mathbf{x}}= \bar{f}(\mathbf{x},\mathbf{u}),
$
where $\bar{f}$ comprises the  last four components in $f$.

In the LQR controllers, we start by defining for each phase P3, P4, and P6, a steady-state 
$
\mathbf{x}_{ref} = ( \beta_{ref}, V_{a_{ref}}, \gamma_{ref}, \theta_{ref}),
$
to linearize around and to set as the reference state. 
We define also  the corresponding steady-state control
$
\mathbf{u}_{ref} = ( F_{p,ref}, w_{q,ref} ).
$

In phases P3 and P6, we set the angles of climb or glide, $\gamma_{climb}$ or $\gamma_{glide}$, respectively.
We also set $\alpha=\alpha_{L}$, $\beta=\beta_{med}$ (an intermediate value along the trajectory of those phases).
In such conditions, we determine $V_a$ such that it satisfies the steady-state conditions 
$
(\dot{V}_{a}, \dot{\gamma}, \dot{\theta} )=(0 ,0, 0)
$.
During phase P4-Loiter, we find the speed $V_a$ satisfying
$
(\dot{\beta}, \dot{V}_{a}, \dot{\gamma}, \dot{\theta} )=(0 ,0, 0, 0) 
$
for 
$
({\beta}, \alpha, {\gamma}, {\theta} )=(\beta_0 ,\alpha_0, 0, 0).
$
Later, in the Simulation Results section, we define the values for the reference states (see Table \ref{table:controller_gains}).

As usual in LQR, we consider the error state $\widetilde{\mathbf{x}}= \mathbf{x} - \mathbf{x}_{ref}$ 
and the error control $\widetilde{\mathbf{u}}= \mathbf{u} - \mathbf{u}_{ref}$.
The linearized model is 
$\dot{\widetilde{\mathbf{x}}} = {A} \widetilde{\mathbf{x}} + {B} \widetilde{\mathbf{u}}$, where
\begin{align}
{A} = \left.\dfrac{d{f} (\widetilde{\mathbf{x}}, \widetilde{\mathbf{u}})}{d\widetilde{\mathbf{x}}}\right|_{\substack{
\widetilde{\mathbf{x}} = 0  \\
\widetilde{\mathbf{u}} = 0   }}, \quad 
{B} = \left.\dfrac{d{f} (\widetilde{\mathbf{x}}, \widetilde{\mathbf{u}})}{d\widetilde{\mathbf{u}}}\right|_{\substack{
\widetilde{\mathbf{x}} = 0  \\
\widetilde{\mathbf{u}} = 0  }} .
\end{align}
In the objective function 
$
\displaystyle \int_0^\infty (\widetilde{\mathbf{x}}^T Q \widetilde{\mathbf{x}} + \widetilde{\mathbf{u}}^T R \widetilde{\mathbf{u}} ) dt
$,
the matrices $Q$ and $R$ are diagonal, initially set using Bryson and Ho \cite{bryson_applied_1975} rule and then manually tuned to achieve the desired response.
In Table IV, we report the matrix values used in the simulation for each phase.
We note that the coefficient in $Q$ associated with the state $\beta$, which varies significantly in P3 and P6, was set to zero in those phases.
\section{Simulation Results}

\subsection{Small-scale Prototype and Simulation Parameters}
\begin{figure}[!ht]
    \centering
    \includegraphics[width=\columnwidth]{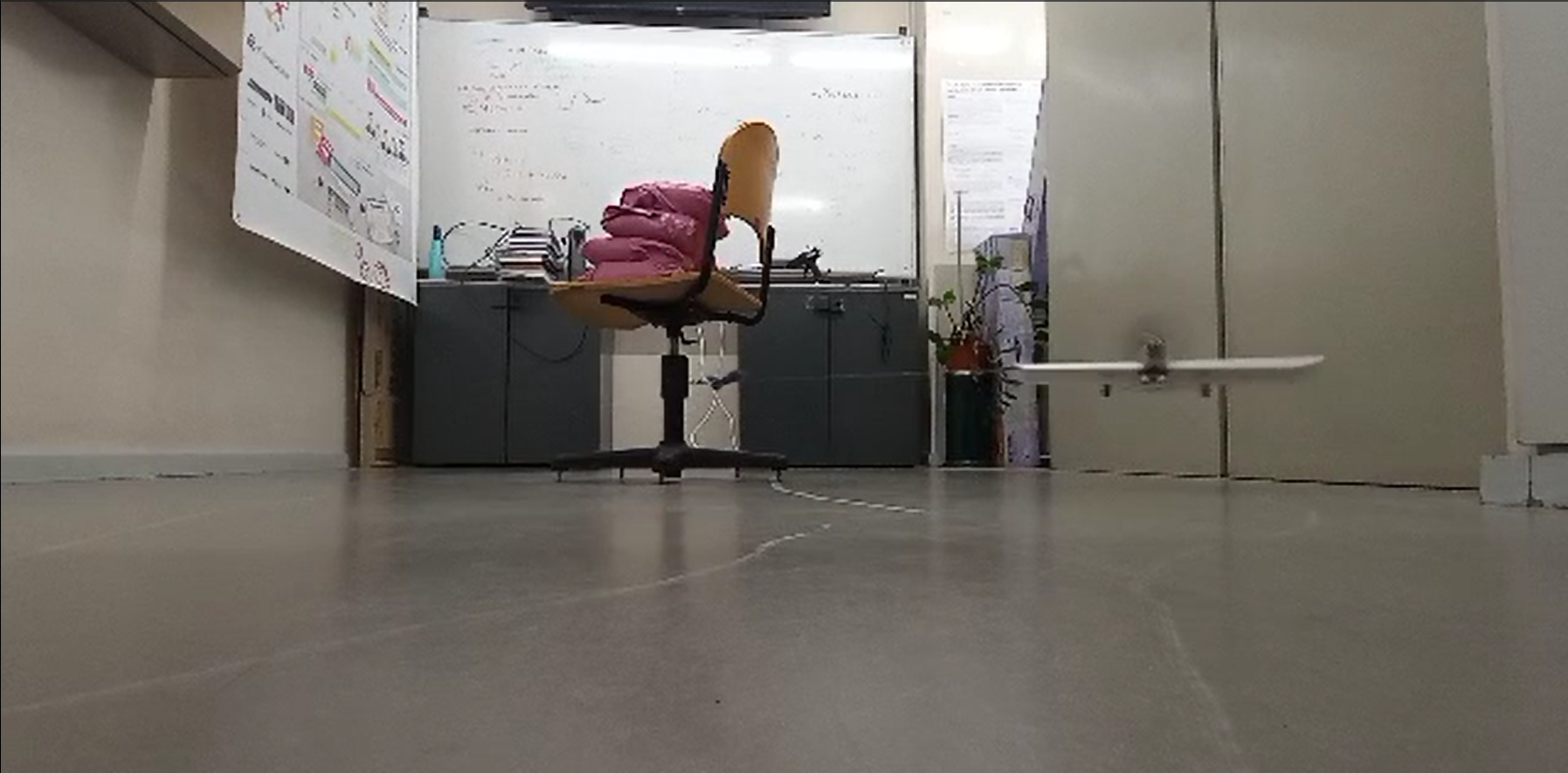}
    \caption{A small-scale prototype during a CTOL flight test validation in our laboratory \cite{fernandes_take-off_2023}.}
    \label{fig:small-scale prototype}
\end{figure}

The simulation parameters are defined taking into account a small-scale aircraft,  to ease the comparison of the simulation results with future experiments (see Fig. \ref{fig:small-scale prototype}). 
The small-scale prototype is equipped with a wing that follows the specifications of NACA 4412 airfoil and is mounted with an incidence angle $\alpha_i= 6^\circ$. 
The maximum $c_{\rm{L}} = 1.4002$ is obtained when: $\alpha_{\text{wing}}=15^\circ$, that corresponds to $\alpha_{\text{aircraft}}=9^\circ$  $(\alpha_{\text{aircraft}}=\alpha_{\text{wing}}-\alpha_i)$.
The maximum $\sfrac{c_{\rm{L}}}{c_{\rm{D}}}= 76.557$ is obtained when: $\alpha_{\text{wing}}=6^\circ$ ($\alpha_{\text{aircraft}}=0^\circ$). Other aircraft and simulation parameters can be seen in Table \ref{tab:parameters}.

\begin{table}[!ht]
\centering

\caption{Aircraft and Simulation Parameters}
\begin{tabular}{c}
\midrule
\textbf{Aircraft and Environment Parameters}\\
\midrule
$A = 0.0576\,\rm{m^2}$, $b = 0.60\,\rm{m}$, $m = 0.350\,\rm{kg}$, $r  = 4b = 2.4\,\rm{m}$ \\
$\omega_q \in \left[-20;20\right]^\circ\,\rm{s^{-1}}$, $F_p \in \left[0;1.5\right]\rm{N}$ \\
$\rho = 1.225\,\rm{kg\,m^{-3}}$, $g = 9.8\,\rm{m\,s^{-2}}$, $V_{w}=0\,\rm{m\,s^{-1}}$\\

\midrule
\textbf{Conditions for phase transition}\\
\midrule
$h_0 = b/2 = 0.3\,\rm{m}$ ($\beta_{0} = 7.18^\circ$) \\
$h_{flare} = 0.063\,\rm{m}$ ($\beta_{flare} = 1.50^\circ$) \\
$V_{rot} = 7.98\,\rm{m\,s^{-1}}$, 
$V_{loiter} = 10.84\,\rm{m\,s^{-1}}$,
$V_{glide} = 8.29\,\rm{m\,s^{-1}}$ \\
$\gamma_{climb} = 3.00^\circ$, 
$\gamma_{glide} = -1.00^\circ$ \\
$\theta_{rot} = 9.00^\circ$,
$\theta_{flare} = 12.00^\circ$ \\

\end{tabular}
\label{tab:parameters}
\end{table}

\begin{table}[!ht]
\caption{Controller references, gains, and objectives for each phase} 
\label{table:controller_gains}
\centering
\begin{tabular}{ll}
\toprule
$\text{\textbf{PID}}_{\text{\textbf{P1}},\theta}$     & $\begin{aligned} (k_p,k_i,k_d) &= (1.00,0.001,0.01) \\ 
\theta_{ref} &= 0^\circ\end{aligned}$ \vspace{5pt}\\

$\text{\textbf{PID}}_{\text{\textbf{P1}},V_a}$     & $\begin{aligned} (k_p,k_i,k_d) &= (0.7,0.08,0.05) \\ 
V_{a_{ref}} &= 7.98\,\rm{m\,s^{-1}}\end{aligned}$ \\
\midrule
\midrule

$\text{\textbf{PID}}_{\text{\textbf{P2}},\theta}$     & $\begin{aligned} (k_p,k_i,k_d) &= (30.00,0.01,1.00) \\ 
\theta_{ref} &= 12.00^\circ\end{aligned}$ \vspace{5pt}\\

$\text{\textbf{PID}}_{\text{\textbf{P2}},V_a}$     & $\begin{aligned} (k_p,k_i,k_d) &= (10.00,0.10,0.01) \\ 
V_{a_{ref}} &= 7.98\,\rm{m\,s^{-1}}\end{aligned}$ \\
\midrule
\midrule

$\text{\textbf{LQR}}_{\text{\textbf{P3}}}$    & $\begin{aligned} &(\beta,V_a,\gamma,\theta)_{_{ref}} = (5.00, 8.25, 3.00, 12.00) \\ 
&Q = \text{diag}([0 \ \ 0.015 \ \ 364.76 \ \ 22.80])\\
&R = \text{diag}([4.83 \ \ 959.18]) \end{aligned}$ \\

\midrule
\midrule

$\text{\textbf{LQR}}_{\text{\textbf{P4}}}$     & $\begin{aligned} &(\beta,V_a,\gamma,\theta)_{_{ref}} = (7.18, 10.84, 0, 0) \\ 
&Q = \text{diag}([0.064 \ \ 0.085e{-3} \ \ 5.62 \ \ 0.033]) \times 10^{3} \\
&R = \text{diag}([2.61 \ \ 8.21]) \end{aligned}$ \\
\midrule
\midrule

$\text{\textbf{PID}}_{\text{\textbf{P5}},\gamma}$     & $\begin{aligned} (k_p,k_i,k_d) &= (9.00, 0.01, 0.10) \\ 
\theta_{ref} &= 9.00^\circ \: \text{and} \: \:\gamma_{ref} = 0^\circ\end{aligned}$\vspace{5pt} \\

$\text{\textbf{PID}}_{\text{\textbf{P5}},V_a}$     & $\begin{aligned} (k_p,k_i,k_d) &= (10.00, 0.10, 1.00) \\ 
V_{a_{ref}} &= 8.29\,\rm{m\,s^{-1}}\end{aligned}$ \\
\midrule
\midrule

$\text{\textbf{LQR}}_{\text{\textbf{P6}}}$     & $\begin{aligned} &(\beta,V_a,\gamma,\theta)_{_{ref}} = (2.39, 7.81, -1.00,  8.00) \\ 
&Q = \text{diag}([0 \ \ 0.015e{-3} \ \ 1.46 \ \ 0.037]) \times 10^{3}\\
&R = \text{diag}([46.91 \ \ 33.29]) \end{aligned}$ \\
\midrule
\midrule

$\text{\textbf{PID}}_{\text{\textbf{P7}},\theta}$     & $\begin{aligned} (k_p,k_i,k_d) &= (1.00,0.01,0.50) \\ 
\theta_{ref} &= 12^\circ (F_p= 0) \end{aligned}$ \\

\bottomrule
\end{tabular}
\end{table}

\subsection{Results}

\begin{figure}[!ht]
    \centering
    \includegraphics[width=\columnwidth,trim={0cm 2.6cm 0cm 2.9cm},clip]{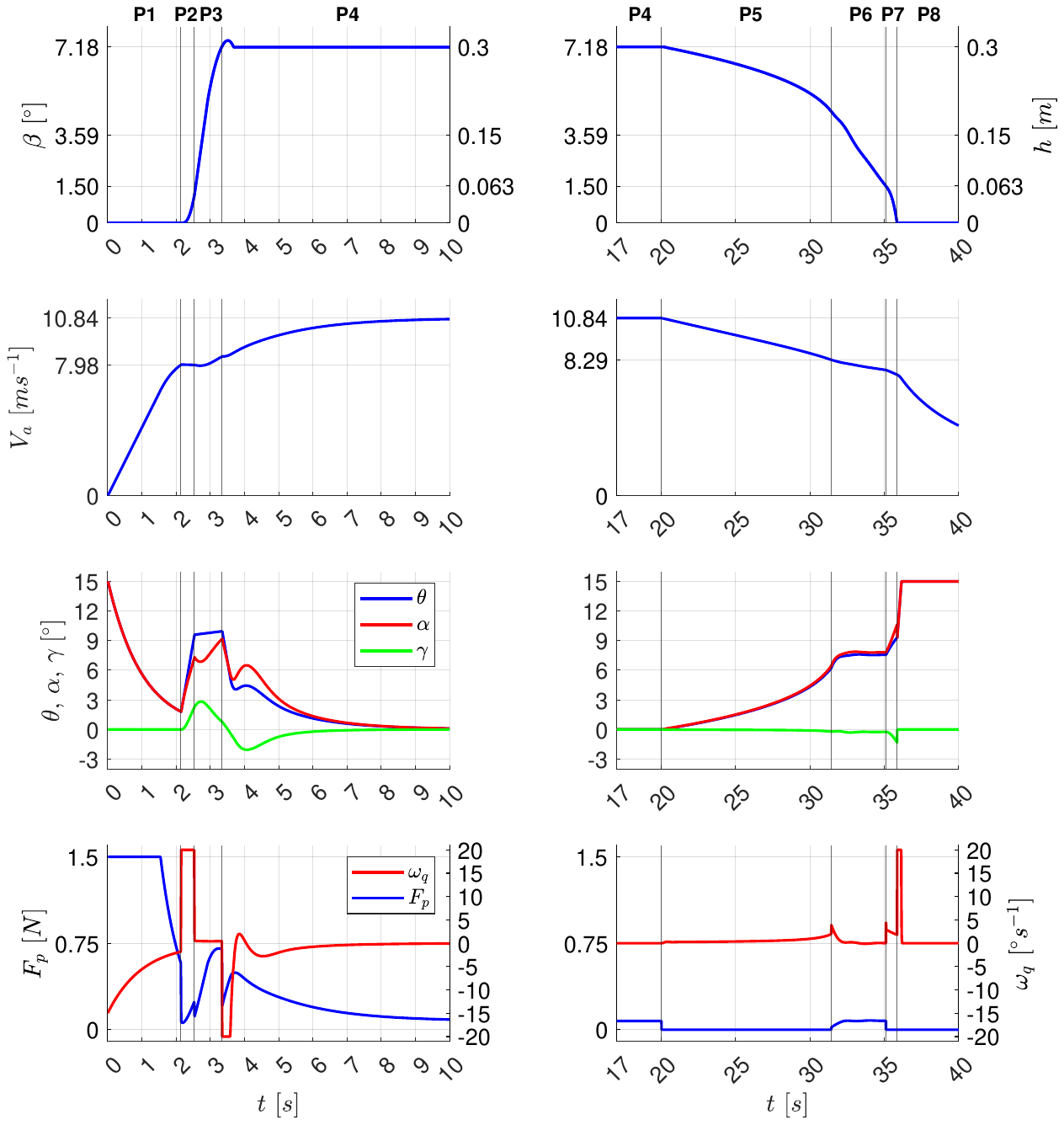}
    \caption{Simulation results during Take-off (left column) and Approach/Landing (right column) phases.  Time intervals for each phase, in seconds:
    $t_{\rm{P1}} \in [0,2.14]$, $t_{\rm{P2}} \in [2.14,2.53]$, $t_{\rm{P3}} \in[2.53,3.34]$, 
     \linebreak$t_{\rm{P4}} \in [3.34,20]$,  $t_{\rm{P5}} \in [20,31.43]$, $t_{\rm{P6}} \in [31.43,35.11]$,
     \linebreak$t_{\rm{P7}} \in [35.11,35.86]$, $t_{\rm{P8}} \geq 35.86$.}
    \label{fig:results}
\end{figure}

A summary of the controllers, their respective reference values, their gains, and objectives, used for each phase is given in Table \ref{table:controller_gains}.
The controllers were implemented using Matlab/Simulink and the results of the simulation are displayed in Fig. \ref{fig:results}.
In the left column, we can see the simulation results during the take-off sub-phases, while the results for the approach and landing sub-phases are given in the right column of the figure. 
The first row of graphs shows the elevation angle of the kite in the left axis and the height in the right axis. 
We can see the speed trajectory in the second row of graphs. In the third row, the attitude angles of the aircraft are plotted. Finally,  the last row of graphs displays the control inputs, with the left axis for the force of the propeller and the right axis for the pitch rate. 

The trajectories follow closely the reference profile given. Starting from the rest position, the tethered aircraft sets its maximum throttle during P1. Then, the pitch rate saturates at its maximum to rapidly tilt up the aircraft while $\theta$ goes to the desired value. The kite starts to gain some altitude during this second phase. After $t=2.53\,\rm{s}$, the aircraft starts the initial climb, while maintaining the desired speed. Finally, after reaching the loiter altitude, we can see a small overshoot. This can be explained by the limitation on the pitch rate that is saturated when trying to reduce the angle-of-attack of the aircraft required for the loiter phase. It then successfully maintains its stabilized altitude ($h_0 = 0.3\,\rm{m}$) during the Loiter phase.
Note that the Loiter phase lasts for almost 18 seconds, and only a few seconds are represented in both graphs. 
The landing command was triggered at $t=20\,\rm{s}$. The approach phase starts by decreasing the speed of the aircraft and setting $F_p$ to zero. After this deceleration, the glide phase starts, maintaining a constant angle-of-attack and constant speed. Finally, a minimum altitude is achieved and the flare maneuver is initiated. Here, the pitch rate jumps to the maximum to rapidly tilt up the aircraft. The aircraft touches the ground and the speed is decreased until it eventually stops. Only a small part of the deceleration to rest is represented in the graphs.

\section{Discussion and Conclusion}

\begin{figure}[!ht]
    \centering
    \includegraphics[width=0.48\textwidth,trim={3.5cm 10cm 4cm 10cm},clip]{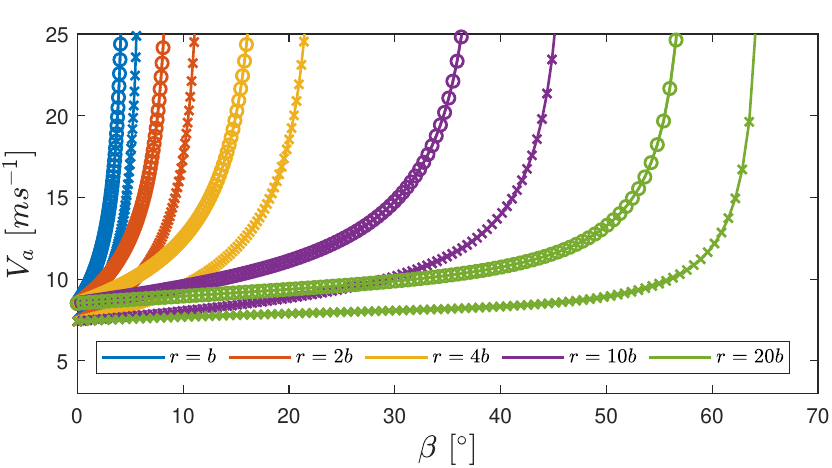}
    \caption{$\beta$ vs $V_a$ during take-off for $\alpha=0$ (curves with circles) and for $\alpha=9^\circ$ (curves with crosses).}
    \label{fig:max_va}
\end{figure}

In general, the results correspond to the expected values, essentially following the reference values previously set.
The aircraft successfully passes through all take-off and landing sub-phases using a simple set of controllers - PIDs and LQRs.
The circular take-off and landing approach was shown to be a viable method to automate such aircraft maneuvers, especially in tethered aircraft as the one used in Airborne Wind Energy Systems.

However, we can point out some aspects that differ from the initial plan. We expected less time to decrease the aircraft speed during phase P5. 
If needed, this can be attenuated by increasing the drag during this phase, for example equipping the aircraft with flaps and deploying them at this stage.

The speed required for a certain elevation angle depends mainly on the aircraft size and on the tether length. 
The drag force and the weight components above a certain speed can be negligible when compared with the lift force and centrifugal force, because they are small. Therefore, taking (\ref{eq:model}) and solving for $\dot{\gamma} = 0$, with ${\gamma} = 0$, we have
for high speeds
%
\begin{equation}
    \label{eq:beta_max}
F_{\rm{L}} = \frac{m {V_a}^2}{r}\cgamma \tbeta \Leftrightarrow
\tbeta = \frac{\frac{1}{2}\rho A c_L }{m}r.
\end{equation}
This equation identifies, for the aircraft dimensions defined previously, the maximum attainable $\beta$ for a specific tether length.

Plotting $V_a$ as a function of the tether length for different values of  $\beta$, using equation $\dot{\gamma} = 0$, we get Fig. \ref{fig:max_va}. From this figure, we can draw some interesting conclusions. We can identify the minimum speed required to start elevating the aircraft for the two values of the AoA considered:  $\alpha = 0^\circ$ and $\alpha=9^\circ$, that correspond to maximum  $\sfrac{c_{\rm{L}}}{c_{\rm{D}}}$ and the maximum lift, respectively.
But, more relevant is that the graph explicitly shows that only a determined $\beta_{max}$ is attainable for certain values of $r$. Also, the high speeds required for this $\beta_{max}$ may not be achievable.  Therefore, the tether length is an important design parameter substantially defining the maximum height that can be obtained using a circular take-off and landing scheme.

\bibliography{ATOL_ECC24}   

\end{document}